\input amstex
\documentstyle{amsppt}
\magnification=\magstep1                        
\hsize6.5truein\vsize8.9truein                  
\NoRunningHeads
\loadeusm

\magnification=\magstep1                        
\hsize6.5truein\vsize8.9truein                  
\NoRunningHeads
\loadeusm

\document
\topmatter

\title
Improved results on the oscillation of the modulus of the Rudin-Shapiro polynomials on the unit circle\\
\endtitle

\rightheadtext{the oscillation of the modulus of the Rudin-Shapiro polynomials on the unit circle}

\author Tam\'as Erd\'elyi
\endauthor

\address Department of Mathematics, Texas A\&M University,
College Station, Texas 77843, College Station, Texas 77843 \endaddress

\endthanks


\date September 12, 2018 
\enddate

\email terdelyi\@math.tamu.edu
\endemail

\dedicatory Dedicated to Paul Nevai on the occasion of his 70th birthday
\enddedicatory

\abstract
Let either $R_k(t) := |P_k(e^{it})|^2$ or $R_k(t) := |Q_k(e^{it})|^2$, where $P_k$ and $Q_k$ are 
the usual Rudin-Shapiro polynomials of degree $n-1$ with $n=2^k$.  
In a recent paper we combined close to sharp upper bounds for the modulus of the autocorrelation coefficients of
the Rudin-Shapiro polynomials with a deep theorem of Littlewood to prove that there is an absolute constant 
$A>0$ such that the equation $R_k(t) = (1+\eta )n$ has at least $An^{0.5394282}$ distinct zeros in 
$[0,2\pi)$ whenever $\eta$ is real, $|\eta| < 2^{-11}$, and $n$ is sufficiently large.  
In this paper we show that the equation $R_k(t)=(1+\eta)n$ has at least $(1/2-|\eta|-\varepsilon)n/2$ 
distinct zeros in $[0,2\pi)$ for every $\eta \in (-1/2,1/2)$, $\varepsilon > 0$, and sufficiently large 
$k \geq k_{\eta,\varepsilon}$. 
\endabstract

\endtopmatter 

\head 1. Introduction \endhead
Let $D := \{z \in {\Bbb C}: |z| < 1\}$ denote the open unit disk of the complex plane.
Let $\partial D :=  \{z \in {\Bbb C}: |z| = 1\}$ denote the unit circle of the complex plane.
The Mahler measure $M_{0}(f)$ is defined for bounded measurable functions $f$ on $\partial D$ by
$$M_{0}(f) := \exp\left(\frac{1}{2\pi} \int_{0}^{2\pi}{\log|f(e^{it})|\,dt} \right)\,.$$
It is well known, see [HL-52], for instance, that
$$M_{0}(f) = \lim_{q \rightarrow 0+}{M_{q}(f)}\,,$$
where
$$M_{q}(f) := \left( \frac{1}{2\pi} \int_{0}^{2\pi}{\left| f(e^{it}) \right|^q\,dt} \right)^{1/q}\,, 
\qquad q > 0\,.$$
It is also well known that for a function $f$ continuous on $\partial D$ we have
$$M_{\infty}(f) := \max_{t \in [0,2\pi]}{|f(e^{it})|} = \lim_{q \rightarrow \infty}{M_{q}(f)}\,.$$
It is a simple consequence of the Jensen formula that
$$M_0(f) = |c| \prod_{j=1}^n{\max\{1,|z_j|\}}$$
for every polynomial of the form
$$f(z) = c\prod_{j=1}^n{(z-z_j)}\,, \qquad c,z_j \in {\Bbb C}\,.$$
See [BE-95, p. 271] or [B-02, p. 3], for instance.

Let ${\Cal P}_n^c$ be the set of all algebraic polynomials of degree at most $n$ with complex coefficients.
Let ${\Cal T}_n$ be the set of all real (that is, real-valued on the real line) trigonometric polynomials
of degree at most $n$. Finding polynomials with suitably restricted coefficients and maximal Mahler measure 
has interested many authors. The classes
$${\Cal L}_n := \left\{ f: \enskip f(z) = \sum_{j=0}^{n}{a_jz^j}\,, \quad a_j \in \{-1,1\} \right\}$$
of Littlewood polynomials and the classes
$${\Cal K}_n := \left\{ f: \enskip f(z) = \sum_{j=0}^{n}{a_jz^j}\,, \quad a_j \in {\Bbb C}, \enskip |a_j| =1 \right\}$$
of unimodular polynomials are two of the most important classes considered.
Observe that ${\Cal L}_n \subset {\Cal K}_n$ and
$$M_0(f) \leq M_2(f) = \sqrt{n+1}\,, \qquad f \in {\Cal K}_n\,.$$
Beller and Newman [BN-73] constructed unimodular polynomials $f_n \in {\Cal K}_n$ such that 
$M_0(f_n) \geq \sqrt{n}-c/\log n$ with an absolute constant $c > 0$.

Section 4 of [B-02] is devoted to the study of Rudin-Shapiro polynomials.
Littlewood asked if there were polynomials $f_{n_k} \in {\Cal L}_{n_k}$ satisfying
$$c_1 \sqrt{n_k+1}  \leq |f_{n_k}(z)| \leq c_2 \sqrt{n_k+1}\,, \qquad z \in \partial D\,,$$
with some absolute constants $c_1 > 0$ and $c_2 > 0$, see [B-02, p. 27] for a reference
to this problem of Littlewood.
To satisfy just the lower bound, by itself, seems very hard, and no such sequence $(f_{n_k})$
of Littlewood polynomials $f_{n_k} \in {\Cal L}_{n_k}$ is known. A sequence of Littlewood polynomials
that satisfies just the upper bound is given by the Rudin-Shapiro polynomials. The Rudin-Shapiro
polynomials appear in Harold Shapiro's 1951 thesis [S-51] at MIT and are sometimes called just
Shapiro polynomials. They also arise independently in Golay's paper [G-51]. They are
remarkably simple to construct and are a rich source of counterexamples to possible
conjectures. The Rudin-Shapiro polynomials are defined recursively as follows:
$$\split P_0(z) & :=1\,, \qquad Q_0(z) := 1\,, \cr 
P_{k+1}(z) & := P_k(z) + z^{2^k}Q_k(z)\,, \cr
Q_{k+1}(z) & := P_k(z) - z^{2^k}Q_k(z)\,, \cr \endsplit$$
for $k=0,1,2,\ldots\,.$ Note that both $P_k$ and $Q_k$ are polynomials of degree $n-1$ with $n := 2^k$
having each of their coefficients in $\{-1,1\}$.
In signal processing, the Rudin-Shapiro polynomials have good autocorrelation properties
and their values on the unit circle are small. Binary sequences with low autocorrelation
coefficients are of interest in radar, sonar, and communication systems.

It is well known and easy to check by using the parallelogram law that
$$|P_{k+1}(z)|^2 + |Q_{k+1}(z)|^2 = 2(|P_k(z)|^2 + |Q_k(z)|^2)\,, \qquad z \in \partial D\,.$$
Hence
$$|P_k(z)|^2 + |Q_k(z)|^2 = 2^{k+1} = 2n\,, \qquad z \in \partial D\,. \tag 1.1$$
It is also well known (see Section 4 of [B-02], for instance), that
$$Q_k(-z) = P_k^*(z) = z^{n-1}P_k(1/z)\,, \qquad z \in \partial D\,,$$
and hence
$$|Q_k(-z)| = |P_k(z)|\,, \qquad z \in \partial D\,. \tag 1.2$$
Let $K := {\Bbb R} \enskip(\text {mod}\,\, 2\pi)$.
Let $m(A)$ denote the one-dimensional Lebesgue measure of $A \subset K$.  
In 1980 Saffari conjectured the following.

\proclaim{Conjecture 1.1}
Let $P_k$ and $Q_k$ be the Rudin-Shapiro polynomials of degree $n-1$ with $n := 2^k$.
We have
$$M_q(P_k) = M_q(Q_k) \sim \frac{2^{k+1)/2}}{(q/2+1)^{1/q}}$$
for all real exponents $q > 0$. Equivalently, we have
$$\split & \lim_{k \rightarrow \infty} 
m{\left(\left\{t \in K: \left| \frac{P_k(e^{it})}{\sqrt{2^{k+1}}} \right|^2 \in [\alpha,\beta] \right\}\right)} \cr 
= \, & \lim_{k \rightarrow \infty}
m{\left(\left\{t \in K: \left| \frac{Q_k(e^{it})}{\sqrt{2^{k+1}}} \right|^2 \in [\alpha,\beta] \right\}\right)} 
= 2\pi(\beta - \alpha) \cr \endsplit$$
whenever $0 \leq \alpha < \beta \leq 1$.
\endproclaim

This conjecture was proved for all even values of $q \leq 52$ by Doche [D-05]
and Doche and Habsieger [DH-04]. Recently B. Rodgers [R-16] proved Saffari's Conjecture 1.1
for all $q > 0$. See also [EZ-17]. An extension of Saffari's conjecture is Montgomery's conjecture below.

\proclaim{Conjecture 1.2}
Let $P_k$ and $Q_k$ be the Rudin-Shapiro polynomials of degree $n-1$ with $n := 2^k$.
We have
$$\split & \lim_{k \rightarrow \infty} 
m{\left(\left\{t \in K: \frac{P_k(e^{it})}{\sqrt{2^{k+1}}} \in E \right\}\right)} \cr
= \, & \lim_{k \rightarrow \infty}
m{\left(\left\{t \in K: \frac{Q_k(e^{it})}{\sqrt{2^{k+1}}} \in E \right\}\right)} 
= 2m(E) \cr \endsplit$$
for any measurable set $E \subset D := \{z \in {\Bbb C}: |z| < 1\}\,.$
\endproclaim

B. Rodgers [R-16] proved Montgomery's Conjecture 1.2 as well.

\head 2. New Results \endhead
Let either $R_k(t) := |P_k(e^{it})|^2$ or $R_k(t) := |Q_k(e^{it})|^2$.
In [AC-18] we combined close to sharp upper bounds for the modulus of the autocorrelation coefficients of 
the Rudin-Shapiro polynomials with a deep theorem of Littlewood (see Theorem 1 in [L-66]) to prove that
there is an absolute constant $A>0$ such that the equation $R_k(t) = (1+\eta)n$ with $n := 2^k$ has at least 
$An^{0.5394282}$ distinct zeros in $K$ whenever $\eta$ is real, $|\eta| \leq  2^{-11}$, and $n$ is sufficiently 
large.  In this paper we improve this result substantially.    

\proclaim{Theorem 2.1}
The equation $R_k(t)=n$ has at least $n/4+1$ distinct zeros in $K$.
Moreover, with the notation $t_j := 2\pi j/n$, there are at least $n/2+2$ values of $j \in \{0,1\ldots,n-1\}$ 
for which the interval $[t_j,t_{j+1}]$ has at least one zero of the equation $R_k(t)=n$.
\endproclaim

\proclaim{Theorem 2.2}
The equation $R_k(t)=(1+\eta)n$ has at least $(1/2-|\eta|-\varepsilon)n/2$ distinct zeros in $K$ 
for every $\eta \in (-1/2,1/2)$, $\varepsilon > 0$, and sufficiently large $k \geq k_{\eta,\varepsilon}$.
\endproclaim

\head 3. Lemma \endhead

In the proof of Theorem 2.1 we need the lemma below stated and proved as Lemma 3.1 in [E-16].   

\proclaim{Lemma 3.1} Let $n \geq 2$ be an integer, $n := 2^k$, and let
$$z_j := e^{it_j}\,, \quad t_j := \frac{2\pi j}{n}\,, \quad j \in {\Bbb Z}\,.$$
We have
$$\split P_k(z_j) & = 2P_{k-2}(z_j)\,, \qquad j=2u\,, \enskip u \in {\Bbb Z}\,, \cr
P_k(z_j) & =(-1)^{(j-1)/2} 2i\,Q_{k-2}(z_j)\,, \qquad j=2u+1\,, \enskip u \in {\Bbb Z}\,, \cr \endsplit$$
where $i$ is the imaginary unit.
\endproclaim

\head 4. Proofs \endhead

\demo{Proof of Theorem 2.1}
Let $k \geq 2$ be an integer. Assume that $R_k(t) = |P_k(e^{it})|^2$, the proof in the case $R_k(t) = |Q_k(e^{it})|^2$ 
is the same. For the sake of brevity let
$$A_j := R_{k-2}(t_j)-n/4\,, \qquad j=0,1,\ldots,n\,.$$  
Using the notation of Lemma 3.1 we study the $(n+1)$-tuple $\langle A_0,A_1,\ldots,A_n \rangle$.   
Observe that $R_{k-2}$ is a real trigonometric polynomial of degree $n/4-1 = 2^k/4-1$, and hence   
$R_{k-2}(t)-n/4$ has at most $n/2-2$ zeros in $K$. Therefore the Intermediate Value Theorem yields that the 
number of sign changes in the $(n+1)$-tuple 
$\langle A_0,A_1,\ldots,A_n \rangle$ is at most $n/2-2$. Thus there are integers 
$$0 \leq j_1 < j_2 < \cdots < j_m \leq n-1$$ 
with $m \geq n - (n/2-2) = n/2 + 2$ such that 
$$A_{j_{\nu}}A_{j_{\nu}+1} \geq 0\,, \qquad \nu = 0,1,\ldots,m\,. \tag 4.1$$   
Using Lemma 3.1 we have either
$$\split 16A_{j_{\nu}}A_{j_{\nu}+1} = & (4(R_{k-2}(t_{j_{\nu}})-n/4))(4(R_{k-2}(t_{j_{\nu}+1})-n/4)) \cr  
= & (4|P_{k-2}(e^{it_{j_{\nu}}})|^2-n)(4|P_{k-2}(e^{it_{j_{\nu}+1}})|^2-n) \cr 
= & (|P_k(e^{it_{j_{\nu}}})|^2-n)(|Q_k(e^{it_{j_{\nu}+1}})|^2-n) \cr 
= & (|P_k(e^{it_{j_{\nu}}})|^2-n)(n-|P_k(e^{it_{j_{\nu}+1}})|^2)\,, \endsplit \tag 4.2$$
or
$$\split 16A_{j_{\nu}}A_{j_{\nu}+1} = & (4(R_{k-2}(t_{j_{\nu}}))-n/4)(4(R_{k-2}(t_{j_{\nu}+1})-n/4)) \cr  
= & (4|P_{k-2}(e^{it_{j_{\nu}}})|^2-n)(4|P_{k-2}(e^{it_{j_{\nu}+1}})|^2-n) \cr 
= & (|Q_k(e^{it_{j_{\nu}}})|^2-n)(|P_k(e^{it_{j_{\nu}+1}})|^2-n) \cr 
= & (|n-|P_k(e^{it_{j_{\nu}}})|^2)(|P_k(e^{it_{j_{\nu}+1}})|^2-n)\,. \endsplit \tag 4.3$$
Combining (4.1), (4.2), and (4.3), we can deduce that
$$(|P_k(e^{it_{j_{\nu}}})|^2-n)(|P_k(e^{it_{j_{\nu}+1}})|^2-n) = -16A_{j_{\nu}}A_{j_{\nu}+1} \leq 0\,, 
\qquad \nu = 0,1,\ldots,m\,.$$ 
Hence the the Intermediate Value Theorem implies that $R_k(t)-n = |P_k(e^{it})|^2-n$ has at least one zero in 
each of the intervals 
$$[t_{j_{\nu}},t_{j_{\nu}+1}]\,, \qquad \nu = 0,1,\ldots,m\,.$$
Recalling that $m \geq n/2 + 2$ we conclude that $R_k(t)-n = |P_k(e^{it})|^2-n$ has at least $m/2 = n/4+1$ 
distinct zeros in $K$.   
\qed \enddemo

\demo{Proof of Theorem 2.2}
This follows from the proof of Theorem 2.1 combined with B. Rodgers's resolution of Saffari's Conjecture 1.1. 
Assume that $R_k(t) = |P_k(e^{it})|^2$, the proof in the case 
$R_k(t) = |Q_k(e^{it})|^2$ is the same. Also, we may assume that $\eta > 0$, the case $\eta = 0$ is contained in 
Theorem 2.1. We use the notation in the proof of Theorem 2.1. Recall that each of the intervals 
$$[t_{j_{\nu}},t_{j_{\nu}+1}]\,, \qquad \nu = 0,1,\ldots,m\,,$$ 
has at least one zero of $R_k$.  On the other hand, by Saffari's Conjecture proved by Rodgers [R-16] we have
$$m(\{t \in K: |R_k(t)-n| \leq |\eta|n\}) < 2\pi (1+\varepsilon)|\eta|$$
for every $\eta \in (-1/2,1/2)$, $\varepsilon > 0$, and sufficiently large $k \geq k_{\eta,\varepsilon}$.
Hence, with the notation
$$B_\eta := \{t \in K: |R_k(t)-n| \leq |\eta|n\}\,,$$
there are at least $m-(1+\varepsilon)|\eta|n$ distinct values of $\nu \in \{1,2,\ldots,m\}$ such that
$$[t_{j_{\nu}},t_{j_{\nu}+1}] \setminus B_\eta \neq \emptyset$$ 
for every $\eta \in (-1/2,1/2)$, and sufficiently large $k \geq k_{\eta,\varepsilon}$. Hence by the Intermediate Value 
Theorem there are at least $m-(1+\varepsilon)|\eta|n$ distinct values of $\nu \in \{1,2,\ldots,m\}$ for which 
$|R_k(t)-n| = |\eta|n$ has a zero in $(t_{j_{\nu}},t_{j_{\nu}+1})$ for every $\eta \in (-1/2,1/2)$, $\varepsilon > 0$, 
and sufficiently large $k \geq k_{\eta,\varepsilon}$. Now observe that (1.1) and (1.2) imply that 
$$|P_k(z)|^2 + |P_k(-z)|^2 = 2n\,, \qquad z \in {\partial D}\,,$$  
that is, 
$$R_k(t)-n = n-R_k(t+\pi)\,, \qquad t \in K\,.$$ 
Hence for every $\eta \in (-1/2,1/2)$ the number of distinct zeros of $|R_k(t)-n| = |\eta|n$ in $K$ is exactly twice 
the number of distinct zeros of $R_k(t) = (1+\eta)n$ in $K$. We conclude that there are at least
$$\frac 12 (m-(1+\varepsilon)\eta n) \geq (1/2-|\eta|-\varepsilon)n/2$$ 
distinct values of $\nu \in \{1,2,\ldots,m\}$ for which $R_k(t)-n = \eta n$ has a zero
in $(t_{j_{\nu}},t_{j_{\nu}+1})$ for every $\eta \in (-1/2,1/2)$, $\varepsilon > 0$, and sufficiently large
$k \geq k_{\eta,\varepsilon}$.
\qed \enddemo

\head 5. Acknowledgement \endhead
The author thanks Stephen Choi for checking the details of this paper carefully before its submission.

\enddocument